\documentclass{IET-Conf-Paper}

\header{
This is an author’s copy supplementing a paper submitted to the 10\textsuperscript{th} Hybrid Power Plants \& Systems Workshop 2026 and copyrighted by IET.
}
\usepackage{wrapfig}
\usepackage{tabularx}
\usepackage{listings}
\usepackage{pgfplots}
\usepackage{tikz}
\usepackage[european,straightvoltages]{circuitikz}
\pgfplotsset{compat=1.18}

\usepackage{booktabs}
\usepackage{mathtools}

\usepackage{bbm}
\begin{document}

\title{EVALUATION OF VARIOUS OBJECTIVE FUNCTIONS FOR OPTIMAL REACTIVE POWER FLOW INCLUDING TRANSFORMER TAP CHANGER OPTIMISATION}

\author{
Gerald Gebhardt\corr, 
Bernd Engel
}

\address{
{elenia Institute for High Voltage Technology and Power Systems, Braunschweig, Germany}
\add{*}{gerald.gebhardt@tu-braunschweig.de}}


\begin{abstract}
Modern distribution grids with high penetration of renewable generation provide substantial flexibility through distributed reactive power sources and transformer tap changers. This high degree of freedom can be exploited for optimisation. 
However, choosing an objective function for optimisation is not trivial, e.~g. minimising grid losses may lead to overvoltages and minimising voltage deviations may lead to higher reactive power flows to neighbouring system operators. 
Thus, this paper deals with the design of an objective function for the centralised optimisation of distributed reactive power sources and transformer tap changers. 
Different objectives for characteristic network quantities are investigated for the optimisation and optimised in a combined manner and separately. The consequences of optimising conflicting target values are then analysed. For the optimisation, various grid usage cases of a 110 kV benchmark power grid from SimBench are examined. The investigated power grid is characterised by a high proportion of renewable energy plants. The optimisation is carried out in a data-driven, object-oriented manner using the interior point method with open source software. At the end of the paper, a meaningful optimisation function with combined weighted objectives is derived and the results are analysed.
\end{abstract}

\maketitle
\section{Introduction}

As part of the German `Energiewende', i.~e. the transition of the energy supply system towards renewable energies, ancillary services must be provided to a greater extent than before by decentralised renewable energy~(DRE)~\cite{doctorthesisRG}. 
This shift affects static voltage maintenance and reactive power flows between neighbouring system operators~(SOs)~\cite{qint}. 
Thus, in power grids with a high share of~DRE, reactive power becomes a crucial control variable for the SO during operation, as demonstrated by the 2025 Iberian Peninsula blackout, during which DRE operated with a fixed power factor~\cite{bo-spain}.
Moreover, when discrete control resources, such as transformer tap changers, are included, finding the optimal operating points requires solving a challenging non-convex mixed-integer non-linear programming (MINLP) problem~\cite{MINLP}.
To determine optimal setpoints for reactive power provisions of DRE, several objectives must be considered. 
However, the optimal operating point regarding one objective typically does not coincide with those for other objectives. 

Therefore, objectives regarding optimal reactive power flow (ORPF) found in other papers are discussed below:
The seminal paper~\cite{opf_og} considers an ORPF for the minimisation of system losses in a small example network consisting of 3 buses. 
To reduce the power losses, the voltage is increased at one bus to the upper bound of 1.2~pu.
The implication of overvoltages regarding system stability are not considered for the example network in~\cite{opf_og}.
In~\cite{opf_intro} active power losses are chosen as objective for ORPF, too.
Voltages maintenance is considered as objective in~\cite{rl-opf}.
In~\cite{rpcacps} the above objectives are termed `financial goals' for real power losses and `technical goals' for voltage maintenance. 
The objectives are added linearly with weighting factors. 
This popular procedure is also employed in~\cite{pu-1}--\cite{pu-end}.
In~\cite{pq-obj} power losses and deviations to reactive power setpoints at each high to medium voltage interface are linearly weighted to form a cost function. 
A sophisticated objective is found in~\cite{pso}:
Reactive power setpoints, 
voltage maintenance and
active power losses of DREs and in the power grid are considered together with a penalty factor for bound violations.
In~\cite{OCoWF} the voltage profile, grid losses, tap position changes and reactive power exchanges are considered. The objectives are summed linearly.
Finally, a quadratic aggregation of objectives is used in~\cite{MB} regarding voltage profiles, equipment loadings, losses 
and setpoints across multiple distribution and transmission SOs (DSOs and TSOs).

Poor performance regarding an omitted objective is not penalisable and weighting the objectives is non-trivial.
This paper addresses this issue with a -- to the best of the author's knowledge -- novel approach of handling different objectives:
An optimisation framework is derived to find a Pareto optimal solution to the MINLP using interdependence analysis to tune the weighting parameters for a specific power grid. 
In this framework objectives are combined quadratically to penalise large deviations and normalised in accordance to the per-unit-system to increase scalability to different network sizes.

To more effectively leverage existing network data and improve solver flexibility, the authors developed their own optimisation framework as follows. 
All network parameters are read directly from \texttt{pandapower}~\cite{pp}. 
The provided optimisers from \texttt{pandapower} are not used, since \texttt{pandapower}'s built-in optimiser suffers from poor convergence behavior~\cite{pp-site} and \texttt{PandaModels.jl} imposes additional overhead by requiring a Julia interface~\cite{pp-site}. 
To avoid these drawbacks, the authors formulate the problem as a nonlinear program in \texttt{Pyomo}~\cite{pyomo} and solve it using \texttt{IPOPT}~\cite{ipopt}, thereby combining the ease of Python data handling with a robust open-source non-linear solver.
To solve the MINLP problem, the transformer tap changers are discretised iteratively with respect to their influence.
The effectiveness of the proposed method is demonstrated on the benchmark network `\texttt{1-HV-mixed--0-no\_sw}' from \texttt{SimBench}~\cite{sb}. 
The developed optimiser is publicly available at~\cite{g-opt}.

The remainder of this paper is organised as follows: 
In Sec.~\ref{sec:pre} the notation is established. 
Sec.~\ref{sec:met} introduces network constraints and the objective function design.
Sec.~\ref{sec:res} displays the results of the methodology for a high voltage benchmark network using 1000 historic network states.
Sec.~\ref{sec:con} concludes the paper and the acknowledgements are found in Sec.~\ref{sec:ack}.

\section{Preliminaries and Notation}\label{sec:pre}
This paper uses the power-invariant positive-sequence, since symmetrical conditions are assumed.
Furthermore, network elements, study cases, objectives and operation points are grouped via sets.
The sets of network elements are $\mathcal B$ for buses, 
$\mathcal E$ for external grids,
$\mathcal G$ for generators, 
$\mathcal L $ for lines, 
$\mathcal M $ for loads 
and $\mathcal T$ for transformers.
Furthermore, 
the set $\mathcal{S}$ represents a variable for network element sets,
the set $\mathcal{C}$ is the set of considered study cases, 
the set $\mathcal{X}$ is the set of operation points
and $\mathcal{O}$ is the set of objectives.
The term `slack' is used as a synonym for the entirety of all external grids. 
The power injected into the slack is
\begin{equation}
 \underline{S}_\text{slack}=\sum_{e\in\mathcal{E}}\underline{S}_e .   
\end{equation}
The signing system provided by \texttt{pandapower}~\cite{pp} is used. Thus, the consumer frame convention (CFC) is applied to $\mathcal{E}$ and $\mathcal{G}$ and the generator frame convention (GFC) is applied to $\mathcal B$, $\mathcal L $ and $\mathcal T$.
Note that the active and reactive powers of generators are displayed with a negative sign when referring to the German standard regarding technical requirements for the high voltage domain~\cite{4120}, since the CFC is used in~\cite{4120}.


\section{Methodology -- Network Constraints and Objective Function Design}\label{sec:met}

\subsection{Network Constraints} \label{sec:met:con}

This self contained subsection shows how data for individual branches and buses is used to formulate network constraints. 
First, the directional $\Pi$-branch model and the equality constraints are established.
Then, the handling of inequality constraints is discussed and a reference apparent power for external grids and a conservative approximation for the serial current between buses in introduced. 

\begin{figure}[ht]
    \centering
    \begin{circuitikz}
    \def\L{1.4cm}
    \def\H{-2.2*\L}
    \draw 
    (0,0) node[above]{$k$} to[short, o-,i=$\underline{I}_{s,ki}$] 
    (\L,0) to[short,i=$\underline{I}_{\text{S},s,ki}$] 
    (2*\L,0) to[generic=$\underline{Y}_{\text{S},s}$]
    (3*\L,0) to[short,-,i<=$\underline{n}_{s}^*\underline{I}_{\text{S},s,ik}$]
    (4*\L,0) to[short,-o,i<=$\underline{n}_{s}^*\underline{I}_{s,ik}$]
    (5*\L,0) node[above]{$i'$};

    \draw 
    (0,\H) to[short,o-o] (5*\L,\H);

    \draw (0,0) to[open, v^>= $\underline{U}_k$] (0,\H);
    \draw (5*\L,0) to[open, v^>= $\frac{\underline{U}_i}{\underline{n}_{s}}$] (5*\L,\H);

    \draw 
    (4*\L,0) to[short,*-,i=$\underline{n}_{s}^*\underline{I}_{\text{P},s,ik}$]
    (4*\L,-\L) to[generic=$\frac{\underline{Y}_{\text{P},s}}2$]
    (4*\L,-2*\L) to[short,-*]
    (4*\L,\H);
    \draw 
    (\L,0) to[short,*-,i=$\underline{I}_{\text{P},s,ki}$]
    (\L,-\L) to[generic=$\frac{\underline{Y}_{\text{P},s}}2$]
    (\L,-2*\L) to[short,-*]
    (\L,\H);
    \end{circuitikz}
    \caption{Directional $\Pi$-branch model of a branch $s\in \mathcal{L} \cup \mathcal T $ defined from bus $k$ to bus $i$.}
    \label{fig:Pi}
\end{figure}
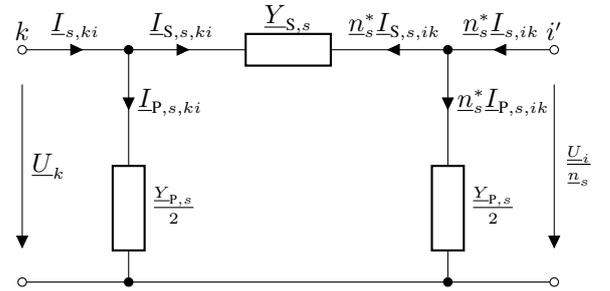

The optimiser discussed in this paper uses the directional $\Pi$-branch model for transformers and lines in accordance with \texttt{pandapower}, since lines have the attributes `\texttt{from\_bus}' and `\texttt{to\_bus}' and transformers have the attributes `\texttt{lv\_bus}' and `\texttt{hv\_bus}'~\cite{pp} there.
Transformers are defined from the lower to the higher voltage side and the orientation of lines is arbitrarily chosen.
The directional $\Pi$-branch is shown in Fig.~\ref{fig:Pi} for a branch $s \in \mathcal{L} \cup \mathcal T$ defined from bus $k$ to bus $i$.
The current $\underline{I}_{s,ki}$ enters $s$ at bus $k$ and is split into a 
current $\underline{I}_{\text{S},s,ki}$ travelling to bus $i$ through $\underline{Y}_{\text{S},s}$ and a 
current $\underline{I}_{\text{S},s,ki}$ travelling to the neutral through $\frac{\underline{Y}_{\text{P},s}}{2}$.
The bus $i$ is transformed via an ideal transformer to the auxiliary bus $i'$ with transmission ratio $\underline{n}_{s}$.
Note that $\underline{n}_{l}=1$ for $l\in \mathcal{L}$ and 
\begin{equation}
    \underline{n}_t=
\frac{U_{\text{N-hv},t}}{U_{\text{N-lv},t}}
    \left[1
    +
    \Delta n_t
     (\psi_t-\psi_{\text{N},t})
     \right]\text{e}^{\text{j}[\varphi_{\text{N},t}
    +
    \Delta\varphi_{n,t}
     (\psi_t-\psi_{\text{N},t})]} \label{e:tap}
\end{equation}
for $t \in \mathcal{T}$, 
with nominal voltages 
$U_{\text{N-hv},t}$ on the higher and 
$U_{\text{N-lv},t}$ on the lower voltage side,
tap step changes regarding magnitude
$\Delta n_t$ and phase $\varphi_{n,t}$ and
present
$\psi_t \in \mathbb{Z}$ and nominal tap position 
$\psi_{\text{N},t}$ as provided by~\cite{pp}.
The optimiser discussed in this paper handles the discreteness of $\psi_t, \forall t\in \mathcal{T}$ as follows:
First an optimisation with $\psi_t \in \mathbbm{R}, \forall t\in \mathcal{T}$ is performed. 
Then, $\psi_t$ for the transformer $t$ with the highest apparent power is rounded to the nearest integer and the value is fixed.
This procedure is repeated until $\psi_t \in \mathbbm{Z}, \forall t\in \mathcal{T}$.
Then a final optimisation is performed. Thus, in total $|\mathcal{T}|+1$ optimisations are executed. 
The admittances 
$\underline{Y}_{\text{S},s}$ and 
$\underline{Y}_{\text{P},s}$ in Fig.~\ref{fig:Pi} are built from \texttt{pandapower}-data.
The entries $\underline{Y}_{ki}$ of the complex bus admittance matrix $\textbf{\textit{Y}}$
are found from these admittances via
\begin{equation}
    \underline{Y}_{ki}= 
    -\sum_{\mathclap{s\in\mathcal{L}\cup\mathcal{T}}}
    \left(
    \frac{\mathbbm{1}_{s,(k,i)}}{\underline{n}_s}
    +
    \frac{\mathbbm{1}_{s,(i,k)}}{\underline{n}_s^*}
    \right)
    \underline{Y}_{\text{S},s}\label{e:Yki}
\end{equation}
for $k\ne i$ and
\begin{equation}
    \underline{Y}_{kk}=
    \sum_{i\in\mathcal{B}\setminus k}
    \sum_{{s\in\mathcal{L}\cup\mathcal{T}}}
    \left(
    \mathbbm{1}_{s,(k,i)}
    +
    \frac{\mathbbm{1}_{s,(i,k)}}{|\underline{n}_s|^2}
    \right)
    \left(
    \underline{Y}_{\text{S},s}
    +\frac{\underline{Y}_{\text{P},s}}{2}
    \right),  \label{e:Ykk}
\end{equation}
with the affiliation operator applied to branches via
\begin{equation}
    \mathbbm{1}_{s,(k,i)} = 
    \left\{\begin{array}{ll}
    1,     &  \text{if $s$ is defined from $k$ to $i$,}\\
    0,     & \text{else.} 
    \end{array}
    \right.
\end{equation}
Equations~(\ref{e:Yki}) and~(\ref{e:Ykk}) can be derived via Kirchhoff's circuit laws as
shown in~\cite{opf_intro} or via the law of conservation of energy as proofed in the appendix.
The entries $\underline{Y}_{ki}$ of the complex bus admittance matrix $\textbf{\textit{Y}}$ are used to formulate the power flow equations 
\begin{equation}
    \sum_{e\in\mathcal{E}}\mathbbm{1}_{e,k} \underline{S}_e
    +\sum_{g\in\mathcal{G}}\mathbbm{1}_{g,k} \underline{S}_g
    =
    \sum_{m\in\mathcal{M}}\mathbbm{1}_{m,k} \underline{S}_m
    +\underline{U}_k\sum_{i\in\mathcal{B}}\underline{Y}_{ki}^*\underline{U}_i^*,
    \label{e:PF:S}
\end{equation}
with the affiliation operator applied to the buses via 
\begin{equation}
    \mathbbm{1}_{s,k} = 
    \left\{\begin{array}{ll}
    1,     &  \text{if $s$ is defined on $k$,}\\
    0,     & \text{else.}
    \end{array}
    \right.
\end{equation}
Furthermore, it is assumed that the voltages at the external buses $\underline{U}_e,\forall e \in \mathcal{E} $ are fixed.

The inequality constraints are 
\begin{align}
   U_{\min,k} \le&\; U_k \le  U_{\text{max},k}, \forall k \in \mathcal{B}, \label{e:con:Uk} \\
  Q_{\min,g}(P_g,U_{k_g})\le&\; Q_g\le Q_{\text{max},g}(P_g,U_{k_g}), \forall g \in \mathcal{G}, \label{e:con:Qg}\\
  -S_{\text{max},e} \le&\; P_e \le S_{\text{max},e}, \forall e \in \mathcal{E},\label{e:con:Pe} \\
  -S_{\text{max},e} \le&\; Q_e \le S_{\text{max},e}, \forall e \in \mathcal{E},\label{e:con:Qe} \\
  0\le & I_{\text{S},s,ki} \le I_{\text{con-max},\text{S},s}, \forall s \in \mathcal{L}\cup\mathcal{T},\label{e:i_conmax,s}\\
   \psi_{\min,t}\le&\;\psi_t \le \psi_{\text{max},t}, \forall t \in \mathcal{T}. \label{e:con:tapt}
\end{align}
As bounds in~(\ref{e:con:Pe}) and~(\ref{e:con:Qe}) the maximal apparent power 
\begin{align}
S_{\text{max},e} =&\, 
\sum_{l \in \mathcal{L}}
\mathbbm{1}_{l,k_e}
U_{\text{max}, l}I_{\text{max}, l}
+
\sum_{t \in \mathcal{T}}
\mathbbm{1}_{t,k_e}
S_{\text{max}, t} \notag
\\
& +
\sum_{m \in \mathcal{M}}
\mathbbm{1}_{m,k_e} S_{\text{N}, m}
+
\sum_{g \in \mathcal{G}}
\mathbbm{1}_{g,k_e} S_{\text{N}, g}
\label{e:S_max,e}
\end{align}
is used. 
To formulate an upper bound for the serial current $\underline{I}_{\text{S},s,ki}$ shown in Fig.~\ref{fig:Pi} in~(\ref{e:i_conmax,s}) the upper bound for the parallel current 
\begin{equation}
I_{\text{max},\text{P},s} = 
    \left|\frac{\underline{Y}_{\text{P},s}}{2}\right|
    \text{max}_{k,i,\underline{n}_{s}} \left(
    U_{\text{max},k}
    ,
    \frac{U_{\text{max},i}}{|\underline{n}_{s}|}
    \right) \label{e:I_max,p}
\end{equation}
is applied via
\begin{equation}
    I_{\text{con-max},\text{S},s} = I_{\text{max},s,ki} -I_{\text{max},\text{P},s}. \label{e:I_max,s}
\end{equation}
This procedure is used to reduce the complexity of the problem, while maintaining feasible operating points.
Equation~(\ref{e:I_max,p}) is suitable,
since usually $\left|\frac{\underline{Y}_{\text{P},s}}{2}\right| \ll |\underline{Y}_{\text{S},s}|$.

\subsection{Objective Function}\label{sec:met:obj}
This subsection presents the objective function design used to formulate the optimisation problem and is organised as follows:
First, the optimisation variables and the general framework for objectives is introduced.
Next, each individual objective is described in detail.
Then, an aggregation procedure for multiple objectives and several study cases is shown.
Finally, a novel evaluation and objective weighting tuning scheme is introduced.

The optimisation variables $\textbf{\textit{x}}$ considered in this paper are $Q_g, \forall g \in \mathcal{G}$ and $\psi_t,\forall t \in \mathcal{T}$, since these variables can be specified by the SO during network operation.
The set of operation point variables is $\mathcal{X}$ and the nomenclature of objectives 
\[
o = (O, \mathcal{S}, A, B)
\]
is as follows:
The Operation $O$ is applied to the set $\mathcal{S}$ to indicate the distance of the objective quantity $A$ to the reference, indicated by the reference indicator $B$.
The objective functions measure the distance of the objective quantity $A_s(\textbf{\textit{x}})$ to the reference $A_{\text{ref},s}$ and scale it by a base quantity $A_{\text{base},s}$ for each member $s$ of the set $\mathcal{S}$. 
The reference indicator $B$ is used to scale the base quantity 
\begin{equation}
    A_{\text{ref},s} = B A_{\text{base},s}.
\end{equation}
The objective functions are applied to certain study cases $c \in \mathcal{C} $ and map to the positive reals including zero, i.~e.
\begin{equation}
    f_o(\textbf{\textit{x}},c) \to \mathbb{R}^+_0,\forall o \in \mathcal{O}.
\end{equation}
Note that a low value for $f_o(\textbf{\textit{x}},c)$ indicates a good performance regarding objective $o$.

The objective functions investigated in this paper are using the operations root-mean-square (rms) and maximum (max) to indicate performance
\begin{align}
    f_{\text{rms},\mathcal{S},A,B}(\textbf{\textit{x}},c) =&\, 
    \sqrt{\frac{1}{|\mathcal{S}|} 
    \sum_{s \in \mathcal{S}} \left(\frac{A_s(\textbf{\textit{x}},c)}{A_{\text{base},s}} -B\right) ^2},
    \\
    f_{\text{max},\mathcal{S},A,B}(\textbf{\textit{x}},c) =&\, \text{max}_{s \in \mathcal{S}}\left| \frac{
    A_s(\textbf{\textit{x}},c)
    }{A_{\text{base}, s}}-B
    \right|.
    \label{e:fmax}
\end{align}

The following objectives are investigated:

1. Maintaining voltage, i.~e.
\begin{equation}
f_{\text{rms},\mathcal{B},U,+1}(\textbf{\textit{x}},c)=
\sqrt{\frac{1}{|\mathcal{B}|} 
    \sum_{k \in \mathcal{B}} \left(\frac{U_k(\textbf{\textit{x}},c)-U_{\text{N},k}}{U_{\text{N},k}}\right) ^2},
\end{equation}
is considered due to the strong coupling between reactive power injections, tap positions and bus voltages.
The base and reference are the nominal voltage $U_{\text{N},k}$ for each bus $k \in \mathcal{B}$.

2. {Minimising reactive power injections of generators}, i.~e.
\begin{equation}
    f_{\text{rms},\mathcal{G},Q,0}(\textbf{\textit{x}},c)=\sqrt{
    \frac{1}{|\mathcal{G}|}\sum_{g\in \mathcal{G}}\left(
    \frac{Q_g(\textbf{\textit{x}},c)}{S_{\text{N},g}}
    \right)^2
    },
\end{equation}
is considered to avoid unnecessary reactive power flows into the power grid. 
However, it is still reactive power injected into the power grid to make the problem feasible, i.~e. to meet the inequality constraints~(\ref{e:con:Uk}) to~(\ref{e:con:tapt}).

3. {Minimising the reactive power's magnitude transferred to the external grids}, i.~e.
\begin{equation}
f_{\text{rms},\mathcal{E},Q,0}(\textbf{\textit{x}},c)=
\sqrt{\frac{1}{|\mathcal{E}|} 
    \sum_{e \in \mathcal{E}} \left(\frac{Q_e(\textbf{\textit{x}},c)}{S_{\text{max},e}}\right) ^2},
    \label{e:rms:Q:E}
\end{equation}
to meet constraints between SOs.
This objective will become more important in the course of Germany's energy transition~\cite{qint}.
The slack is not considered in~(\ref{e:rms:Q:E}), since the ideal of local reactive power compensation, i.~e. $Q_e=0 ,\forall e \in \mathcal{E}$, is assumed.
Note that with this objective transformer loadings of transformers directly connected to the external grid are also considered, since reactive current is minimised. 
The reference can be used as setpoint if reactive power provision to neighbouring SOs is desired.

4. {Minimising active power losses in the power grid}, i.~e.
\begin{equation}
    f_{\text{rms},\text{slack},P,-1}(\textbf{\textit{x}},c)=\left|
    \frac{P_\text{slack}(\textbf{\textit{x}},c)+\sum_{e\in\mathcal{E}}S_{\text{max},e}}{\sum_{e\in\mathcal{E}}S_{\text{max},e}}
    \right|,
    \label{e:f:slackP}
\end{equation}
is considered for economic reasons.
For this objective, the active power entering the slack $P_\text{slack}$ is minimised~\cite{opf_og}. 
The base of the slack is provided by a summation of~(\ref{e:S_max,e}), $\forall e \in \mathcal{E}$.
The objective is not directly applied to individual external grids, e.~g. via $f_{\text{rms},\mathcal{E},P,-1}(\textbf{\textit{x}},c)$, since this would scale the powers differently for each individual external grid.
Note that the rms-operation is equivalent to an evaluation of the magnitude, since the cardinality of the slack is one.

5. {Minimising the maximal loading for
lines}, i.~e.
\begin{equation}
    f_{\text{max},\mathcal{L},I_\text{S},0}(\textbf{\textit{x}},c) = \text{max}_{l \in \mathcal{L}}\left( \frac{
    I_{\text{S}, l}(\textbf{\textit{x}},c)}{I_{\text{con-max,S}, l}}
    \right),
\end{equation}
is considered to penalise operations with high line loadings.
Note that the objective is applied to the serial current to simplify the problem and thus~(\ref{e:i_conmax,s}) is used as base. 

There are many other objectives worth to consider. 
However, this paper focuses on the objectives above, since with this procedure, controllable active and reactive powers, voltages and currents are considered for all operating resources in the power grid. 

In the following an aggregation scheme for the objectives above is presented.
Different objectives $f_o(\textbf{\textit{x}},c)$ are composed with weighting parameters $\alpha_o$ for $o \in \mathcal{O}$
to formulate the main objective 
\begin{equation}
f_\alpha(\textbf{\textit{x}},c) =
\sqrt{
\frac{\sum_{o \in \mathcal{O}} \alpha_o f_o^2(\textbf{\textit{x}},c)}{\sum_{o \in \mathcal{O}} \alpha_{o}}} \label{e:f_a}
\end{equation}
for the optimiser as Pareto optimal solution of multiple objectives.
Note that if only the objective $o$ is considered $f_\alpha(\textbf{\textit{x}},c)$ simplifies to $f_o(\textbf{\textit{x}},c)$.
In~(\ref{e:f_a}) a quadratic composition of objective functions $f_o(\textbf{\textit{x}},c)$ has been chosen, since it gives greater weight to objectives with a larger deviations. 
The optimiser provides the optimal operating point
\begin{equation}
\textbf{\textit{x}}_\alpha^*(c)=\arg \left(
\min_{\textbf{\textit{x}}}
    f_\alpha(\textbf{\textit{x}},c)\right)
    \label{e:xaopt}
\end{equation}
depending on the current weighting $\alpha$ and study case $c$.
Furthermore, the optimal operation point $\textbf{\textit{x}}^*_{\alpha'}$ regarding another weighting ${\alpha'}$ can be evaluated for the weighting $\alpha$ via $f_{\alpha}(\textbf{\textit{x}}^*_{\alpha'},c)$.
Similarly, to evaluate several optimal study points $\mathcal{X}_{{\alpha'}}^*$ obtained for a weighting ${\alpha'}$ regarding another weighting $\alpha$ the performance measure
\begin{equation}
F_{\alpha}(\mathcal{X}_{{\alpha'}}^*,\mathcal{C})=
\sqrt{
\frac{1}{|\mathcal{C}|} 
\sum_{c\in \mathcal{C}}f^2_\alpha(x_{{\alpha'}}^*(c),c) }\label{e:Fab}
\end{equation}
is introduced.
A quadratic composition has been chosen in~(\ref{e:Fab}), since this aligns with~(\ref{e:f_a}).
In~(\ref{e:Fab}) a scaling by the cardinality $|\mathcal{C}|$ is used to normalise the objective in compliance with the per-unit-system.
Note that if only the objective $o$ is considered in the evaluation, (\ref{e:Fab}) simplifies to
\begin{equation}
F_o(\mathcal{X}_{{\alpha'}}^*,\mathcal{C})=
\sqrt{
\frac{1}{|\mathcal{C}|} 
\sum_{c\in \mathcal{C}}f^2_o(x_{{\alpha'}}^*(c),c) }.
\label{e:Fob}
\end{equation}
Furthermore, (\ref{e:Fab}) can be restored from~(\ref{e:Fob}), $\forall o \in \mathcal{O}$ via
\begin{equation}
F_\alpha(\mathcal{X}_{{\alpha'}}^*,\mathcal{C}) 
=
\sqrt{
\frac{\sum_{o \in \mathcal{O}} \alpha_o F_o^2(\mathcal{X}_{{\alpha'}}^*,\mathcal{C})}{\sum_{o \in \mathcal{O}} \alpha_{o}}}.
\end{equation}
The optimal operation points for $\mathcal{C}$ are found from
\begin{equation}
    \mathcal{X}_\alpha^*(\mathcal{C})=\arg \left(
\min_\mathcal{X}
    F_\alpha(\mathcal{X},\mathcal{C})\right)
\end{equation}
similar to~(\ref{e:xaopt}).
The performance of the optimal operation points can be evaluated for all objectives separately via~(\ref{e:Fob}).

\begin{figure}[ht]
\small 
    \centering
    \begin{tikzpicture}
    \def\l{3.5cm}
    \def\r{3.5cm}
    \def\h{1cm}
      \node (opf) [
        draw, 
        rectangle,
        inner sep=0pt, 
        text width = \l, 
        align=center,
        minimum height=\h
      ] at (0,0) {$\min_{\mathcal{X}} F_{o}(\mathcal{X},\mathcal{C}),\break\forall o \in \mathcal{O}$};
      
      \node (pf) [
        draw, 
        rectangle,
        inner sep=0pt, 
        text width = \l, 
        align=center,
        minimum height=\h
      ] at (0,-2*\h) {$F_{o'}(\mathcal{X}_o^*(\mathcal{C}),\mathcal{C}),\break \forall o\in \mathcal{O},o' \in \mathcal{O}$};
      
      \node (opf_a) [
        draw, 
        rectangle,
        inner sep=0pt, 
        text width = \r, 
        align=center,
        minimum height=\h
      ] at (4.5,0) {$\min_{\mathcal{X}} F_\alpha(\mathcal{X},\mathcal{C})$};
      
      \node (pf_a) [
        draw, 
        rectangle,
        inner sep=0pt, 
        text width = \r, 
        align=center,
        minimum height=\h
      ] at (4.5,-2*\h) {$F_{o}(\mathcal{X}_\alpha^*(\mathcal{C}),\mathcal{C}),\break\forall o\in \mathcal{O}$};
    
      \draw[->] (opf) -- node[right]{$\mathcal{X}_{{o}}^*, \forall o \in \mathcal{O}$} (pf);
      \draw[->] (pf.east) -- (2.25,-2*\h) -- node[right, text width = 2cm, align = left] {$\alpha$ from \break analysis} (2.25,0) --  (opf_a.west);
      \draw[->] (opf_a) -- node[right]{$\mathcal{X}^*_\alpha$} (4.5,-1.5*\h);
    \end{tikzpicture}
    \caption{The evaluation pipeline for a set of study cases $\mathcal{C}$, objective function and objective set $\mathcal{O}$.}
    \label{fig:pipeline}
\end{figure}
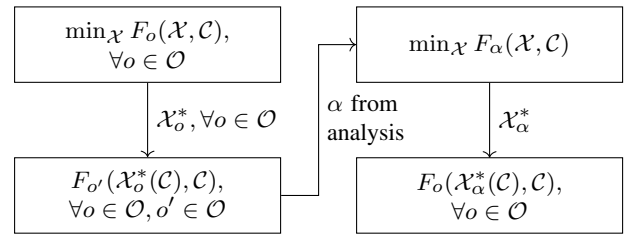

For this optimisation, the SO has to tune the weighting parameters $\alpha$ to align with a specific grid.
To make tuning the parameters $\alpha$ less arbitrary, the procedure shown in Fig.~\ref{fig:pipeline} is proposed.
The first row of Fig.~\ref{fig:pipeline} displays the generation of optimal operation points, whereas the second row shows their evaluation.
In the first column of Fig.~\ref{fig:pipeline}, all objectives are dealt with separately to evaluate the interdependence of different objectives. 
The results are then used to tune the weighting parameters $\alpha$ to cumulate the objectives in second column.
In this work, the expected values $\mu_{F_o}, \forall o \in \mathcal{O}$ are used as marker for average performance to scale the weighting parameters via
\begin{equation}
    \alpha_o = \frac{\tilde{\alpha}_o}{\mu_{F_o}}.
    \label{e:atilde}
\end{equation}
Thus, the SO can tune $\tilde\alpha$ on normalised objectives to align with a specific grid.

\section{Results -- Objective Function Design of a High Voltage Network}
\label{sec:res}
In this section, the evaluation scheme presented in Fig.~\ref{fig:pipeline} is applied to the power grid `\texttt{1-HV-mixed--0-no\_sw}' provided by \texttt{SimBench}~\cite{sb}.
This open source benchmark grid is selected, since it offers a high proportion of DRE and data for $P_g, \forall g \in \mathcal{G}$,
$P_m,Q_m, \forall m \in \mathcal{M}$ and
$U_e, \forall e \in \mathcal{E}$
is provided for one year with a spacing of 15~Minutes between data points.
The network is reduced within the optimiser to accelerate the computations: Admittances of parallel branches are added and network elements with a nominal apparent power of zero are omitted.
The optimiser~\cite{g-opt} uses the interior point method via \texttt{IPOPT}~\cite{ipopt} through the \texttt{Pyomo}-interface~\cite{pyomo} on \texttt{SimBench}-data~\cite{sb} accessed via \texttt{pandapower}~\cite{pp}.

To improve performance, a power flow is performed to initialise the starting values for the optimisation.
The set of 1000 study points $\mathcal{C}$ is constructed using the random seed 42 in \texttt{NumPy} and the function \texttt{random.choice} to avoid duplicates~\cite{np}.
The inequality constraints~(\ref{e:con:Uk}),~(\ref{e:con:Qg}) and~(\ref{e:con:tapt}) are specified as
\begin{align}
   0.9 U_{\text{N},k} \le&\; U_k \le  1.1 U_{\text{N},k}, \forall k \in \mathcal{B}, \\
  -0.5P_{\text{max},g}\le&\; Q_g\le 0.5P_{\text{max},g}, \forall g \in \mathcal{G} ,\label{e:in:Qg}\\
   -16\le&\;\psi_t \le 16, \forall t \in \mathcal{T}.
\end{align}
Note that the reactive power limits retrievable by the DSO for 110~kV grids in Germany depend on the current active power supply and voltage at the connection point~\cite{4120}.
These constraints are not used in (\ref{e:in:Qg}) to make the results more general, to evaluate operating points outside the standard and to improve feasibility.
Especially the STATCOM-domain, i.~e. 
$P_g \le 0.1P_{\text{max},g}$, is of particular interest, since the retrievability of reactive power is handled differently across different version of~\cite{4120}.
Note that the proposed operating points remain feasible if they are technically realisable and the plant operators participate in Germany's reactive power market~\cite{qmarket}.

Since the transformers are directly connected to the external grids in `\texttt{1-HV-mixed--0-no\_sw}',~(\ref{e:S_max,e}) simplifies to
\begin{align}
S_{\text{max},e} =&\, 
\sum_{t \in \mathcal{T}}
\mathbbm{1}_{t,k_e}
S_{\text{max}, t} 
\label{e:S_max,e,app}.
\end{align}

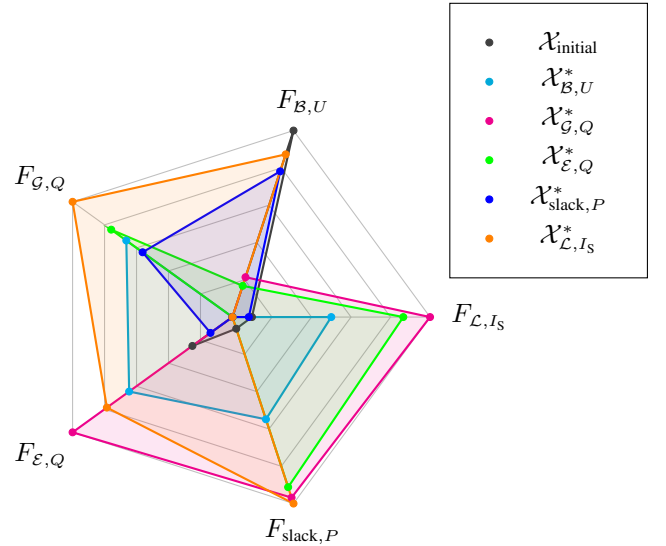
\begin{figure}[ht]
    \centering
\begin{tikzpicture}[scale=2.6]
  \def\n{5}
  \pgfmathsetmacro\angle{360/\n}

  \def\labels{
  {$F_{\mathcal{B},U}$},
  {$F_{\mathcal{G},Q}$},
  {$F_{\mathcal{E},Q}$},
  {$F_{\text{slack},P}$},
  {$F_{\mathcal{L},I_\text{S}}$}
  }
 
\foreach \r in {0.2,0.4,0.6,0.8,1.0} {
    \draw[gray!50]
    ({0*\angle}:\r)
    \foreach \i in {1,...,\n} {
        -- ({\i*\angle}:\r)
    }
    -- cycle;
}

  \foreach [count=\i] \lab in \labels {
    \draw[gray!50] (0,0) -- ({\i*\angle}:1);
    
    \pgfmathparse{int(\i==2||\i==3)}
    \ifnum\pgfmathresult=1
    \node at ({\i*\angle}:1.2) {\lab};
    \fi
    
    \pgfmathparse{int(\i==1||\i==4)}
    \ifnum\pgfmathresult=1
    \node at ({\i*\angle}:1.15) {\lab};
    \fi

    \ifnum\i=5
    \node at ({\i*\angle}:1.25) {\lab};
    \fi
  }

\draw[black](1.1,1.6) -- (2.1,1.6) -- (2.1,0.2) -- (1.1,0.2) --  cycle;

\draw[black!50!gray, thick, fill=black!50!gray, fill opacity=0.1]
(72.0:1.0) -- 
(144.0:0.0) -- 
(216.0:0.2505379322057163) -- 
(288.0:0.062413988574332176) -- 
(360.0:0.09869028959305773) -- 
cycle;
\filldraw[black!50!gray] (1.3,1.4) circle (0.5pt);
\node at (1.7000000000000002,1.4) {$\mathcal{X}_\text{initial}$};
\draw[cyan!75!gray, thick, fill=cyan!75!gray, fill opacity=0.1]
(72.0:0.0) -- 
(144.0:0.6638916336729618) -- 
(216.0:0.6461682600329957) -- 
(288.0:0.548182125435131) -- 
(360.0:0.5000425661177801) -- 
cycle;
\filldraw[cyan!75!gray] (1.3,1.2) circle (0.5pt);
\node at (1.7000000000000002,1.2) {{$\mathcal{X}^*_{\mathcal{B}, U}$}};
\draw[magenta, thick, fill=magenta, fill opacity=0.1]
(72.0:0.21387191300237485) -- 
(144.0:0.0005647647250726469) -- 
(216.0:1.0) -- 
(288.0:0.9681690076776472) -- 
(360.0:1.0) -- 
cycle;
\filldraw[magenta] (1.3,1.0) circle (0.5pt);
\node at (1.7000000000000002,1.0) {{$\mathcal{X}^*_{\mathcal{G}, Q}$}};
\draw[green, thick, fill=green, fill opacity=0.1]
(72.0:0.16744289894085818) -- 
(144.0:0.7583772112591309) -- 
(216.0:0.0) -- 
(288.0:0.9120190193197042) -- 
(360.0:0.8646292595581909) -- 
cycle;
\filldraw[green] (1.3,0.8) circle (0.5pt);
\node at (1.7000000000000002,0.8) {{$\mathcal{X}^*_{\mathcal{E}, Q}$}};
\draw[blue, thick, fill=blue, fill opacity=0.1]
(72.0:0.7813640316003357) -- 
(144.0:0.5620968350176109) -- 
(216.0:0.13737020493399488) -- 
(288.0:0.0) -- 
(360.0:0.08279837796385481) -- 
cycle;
\filldraw[blue] (1.3,0.6) circle (0.5pt);
\node at (1.7000000000000002,0.6) {{$\mathcal{X}^*_{\text{slack}, P}$}};
\draw[orange, thick, fill=orange, fill opacity=0.1]
(72.0:0.8724195076067045) -- 
(144.0:1.0) -- 
(216.0:0.7854561086168497) -- 
(288.0:1.0) -- 
(360.0:0.0) -- 
cycle;
\filldraw[orange] (1.3,0.4) circle (0.5pt);
\node at (1.7000000000000002,0.4) {{$\mathcal{X}^*_{\mathcal{L}, I_\text{S}}$}};

\filldraw[black!50!gray] ({72.0}:1.0) circle (0.5pt);
\filldraw[black!50!gray] ({144.0}:0.0) circle (0.5pt);
\filldraw[black!50!gray] ({216.0}:0.2505379322057163) circle (0.5pt);
\filldraw[black!50!gray] ({288.0}:0.062413988574332176) circle (0.5pt);
\filldraw[black!50!gray] ({360.0}:0.09869028959305773) circle (0.5pt);

\filldraw[cyan] ({72.0}:0.0) circle (0.5pt);
\filldraw[cyan] ({144.0}:0.6638916336729618) circle (0.5pt);
\filldraw[cyan] ({216.0}:0.6461682600329957) circle (0.5pt);
\filldraw[cyan] ({288.0}:0.548182125435131) circle (0.5pt);
\filldraw[cyan] ({360.0}:0.5000425661177801) circle (0.5pt);

\filldraw[magenta] ({72.0}:0.21387191300237485) circle (0.5pt);
\filldraw[magenta] ({144.0}:0.0005647647250726469) circle (0.5pt);
\filldraw[magenta] ({216.0}:1.0) circle (0.5pt);
\filldraw[magenta] ({288.0}:0.9681690076776472) circle (0.5pt);
\filldraw[magenta] ({360.0}:1.0) circle (0.5pt);

\filldraw[green] ({72.0}:0.16744289894085818) circle (0.5pt);
\filldraw[green] ({144.0}:0.7583772112591309) circle (0.5pt);
\filldraw[green] ({216.0}:0.0) circle (0.5pt);
\filldraw[green] ({288.0}:0.9120190193197042) circle (0.5pt);
\filldraw[green] ({360.0}:0.8646292595581909) circle (0.5pt);

\filldraw[blue] ({72.0}:0.7813640316003357) circle (0.5pt);
\filldraw[blue] ({144.0}:0.5620968350176109) circle (0.5pt);
\filldraw[blue] ({216.0}:0.13737020493399488) circle (0.5pt);
\filldraw[blue] ({288.0}:0.0) circle (0.5pt);
\filldraw[blue] ({360.0}:0.08279837796385481) circle (0.5pt);

\filldraw[orange] ({72.0}:0.8724195076067045) circle (0.5pt);
\filldraw[orange] ({144.0}:1.0) circle (0.5pt);
\filldraw[orange] ({216.0}:0.7854561086168497) circle (0.5pt);
\filldraw[orange] ({288.0}:1.0) circle (0.5pt);
\filldraw[orange] ({360.0}:0.0) circle (0.5pt);

\end{tikzpicture}
    \caption{Normalised evaluation of the interdependence of different objectives $F_o$ for operating points $\mathcal{X}_o$.}
    \label{fig:norm:obj}
\end{figure}

In the following, the objectives discussed in subsection~\ref{sec:met:obj} are evaluated. 
For reasons of space, the objectives $o$ are abbreviated by $(\mathcal{S},A)$ for the remainder of this document. 
Furthermore, the pseudo-objective `initial' is introduced and refers to the condition with $Q_g=0,\forall g\in\mathcal{G}$ and $\psi_t=\psi_{\text{N},t},\forall t\in\mathcal{T}$ as baseline for the results. 
Note that for the pseudo-objective `initial' some points are infeasible regarding the inequality constraints, e.~g. due to voltage violations. 
This distinguishes the operation points $\mathcal{X}^*_{\mathcal{G},Q}$ from $\mathcal{X}^*_\text{initial}$, since the operation points $\mathcal{X}^*_{\mathcal{G},Q}$ are feasible.
Table~\ref{tab:res} shows the results regarding the first column of Fig.~\ref{fig:pipeline}. 
\begin{table}[h]
\caption{Evaluation of the interdependence of different objectives $F_o$ for operating points $\mathcal{X}_o$. 
}
\centering
\begin{tabularx}{\linewidth}{XXXXXX}
\toprule
$F_o$
& $F_{\mathcal{B},U}$
& $F_{\mathcal{G},Q}$
& $F_{\mathcal{E},Q}$
& $F_{\text{slack},P}$
& $F_{\mathcal{L},I_\text{S}}$
\\
\midrule
$\mathcal{X}_\text{initial} $
& \textbf{0.08985} & $\underline{\text{0.00000}}$ & 0.03345 & 0.87591 & 0.30937 \\
$\mathcal{X}^*_{\mathcal{B},U}$ 
& $\underline{\text{0.01837}}$ & 0.14906 & 0.08627 & 0.87671 & 0.35196 \\
$\mathcal{X}^*_{\mathcal{G},Q}$ 
& 0.03366 & 0.00013 & \textbf{0.13351} & 0.87740 & \textbf{0.40500} \\
$\mathcal{X}^*_{\mathcal{E},Q}$ 
& 0.03034 & 0.17027 & $\underline{\text{0.00000}}$ & 0.87731 & 0.39064 \\
$\mathcal{X}^*_{\text{slack},P}$  
& 0.07422 & 0.12620 & 0.01834 & $\underline{\text{0.87581}}$ & 0.30769 \\
$\mathcal{X}^*_{\mathcal{L},I_\text{S}}$ 
& 0.08073 & \textbf{0.22452} & 0.10487 & \textbf{0.87745} & $\underline{\text{0.29890}}$ \\
\midrule
$\mu_{F_o}$ 
& 0.05453 & 0.11170 & 0.06274 & 0.87677 & 0.34393 \\
$\sigma_{F_o}$ 
& 0.03050 & 0.09240 & 0.05311 & 0.00075 & 0.04585
\\
\bottomrule
\end{tabularx}
\label{tab:res}
\end{table}
The columns display the performance regarding the objectives and
the optimisation objectives are shown in the first six rows.
The last two rows are statistical analysis regarding the performance. 
The worst performance is bold and the best performance is underlined.
Fig.~\ref{fig:norm:obj} displays the results from Tab.~\ref{tab:res} normalised by the smallest and largest entry for each criteria to improve readability as radar chart. 
A good performance is characterised by a small area in Fig.~\ref{fig:norm:obj}.
In the following, the results are analysed column wise for Table~\ref{tab:res} and counter clockwise for Fig.~\ref{fig:norm:obj}.

The worst performance regarding the objective `maintaining voltage $(\mathcal{B},U)$' is obtained for the `initial'-operation points $\mathcal{X}_{\text{initial}}$, since infeasible points regarding voltage violations are included. 
A bad performance is also attained for $\mathcal{X}^*_{\text{slack},P}$ and $\mathcal{X}^*_{\mathcal{L},I_\text{S}}$, since over-voltage reduces network losses and loadings.
The best performance is achieved for $\mathcal{X}^*_{\mathcal{B},U}$ with a rms-deviation of 1.873\%. 
Excluding deviations at the external net, the rms-deviation for $\mathcal{X}^*_{\mathcal{B},U}$ is 0.089\%.

Regarding the objective `minimising reactive power feed-ins from generators $(\mathcal{G},Q)$' the operation points $\mathcal{X}_{\text{initial}}$ perform best, since the objective is by design zero, even if the operation point is infeasible due to boundary violations.
The second best results are attained for $\mathcal{X}^*_{\mathcal{G},Q}$, since reactive power is only injected to make $c$ feasible.
The operation points $\mathcal{X}^*_{\mathcal{L},I_\text{S}}$ perform worst. 
Thus reducing maximal line's loading results in large amounts of reactive power injections into the considered grid.

The best performance regarding the objective `minimising reactive power feed-ins into the external grids $(\mathcal{E},Q)$' is achieved for the operation points $\mathcal{X}^*_{\mathcal{E},Q}$.
Interestingly, the operation points $\mathcal{X}^*_{\mathcal{G},Q}$ show the worst performance with respect to this objective.

The objective `minimising active power losses in the power grid $(\text{slack},P)$' shows only a small standard deviation\break $\sigma_{F_0}=7.5\cdot 10^{-4}$, i.~e. 1.5~MW, in the results. 
The operation points $\mathcal{X}^*_{\mathcal{L},I_\text{S}}$ show the worst performance regarding this objective.

The small performance increase of $\mathcal{X}_\mathcal{L,I_\text{S}}$ regarding objective `minimising the maximal loading of lines $(\mathcal{L},I_\text{S})$' comes with huge costs.
Thus, it is not a reasonable objective for the considered grid.

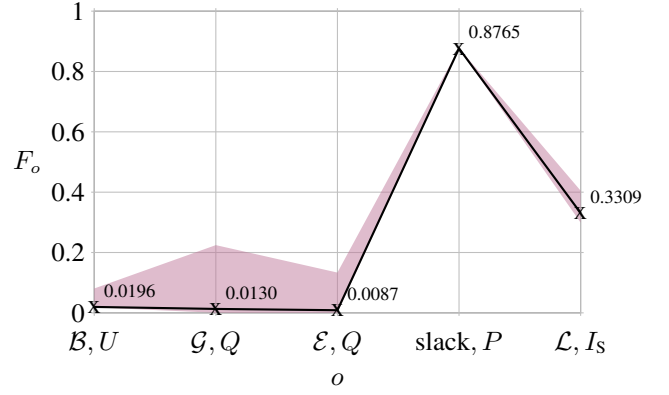
\begin{figure}[ht]
    \centering
\begin{tikzpicture}[scale=4]
  \def\n{5}
  \pgfmathsetmacro\angle{360/\n}

  \def\xlabels{
  {${\mathcal{B},U}$},
  {${\mathcal{G},Q}$},
  {${\mathcal{E},Q}$},
  {${\text{slack},P}$},
  {${\mathcal{L},I_\text{S}}$}
  }

 \def\dx{0.4}
 \def\xshift{0}
  \draw[gray] (0,0) rectangle (\dx*4 + \xshift*2,1) ;

  \foreach \i in {0, 0.2, 0.4, 0.6, 0.8, 1} {
    \draw[gray!50] (-0.01, \i) -- (\dx*4 + \xshift*2 + 0.01, \i);
    \node[left] at (0, \i) {\i};
  }

  \node at (2*\dx + \xshift, -0.22) {$o$};
  \node at (-0.22, 0.5) {$F_o$};

  \foreach [count=\i from 0] \lab in \xlabels {
    \draw[gray!50] (\i*\dx + \xshift,-0.01) -- (\i*\dx + \xshift,1.01);
    \node at (\i*\dx + \xshift,-0.1) {\lab};
  }

  \fill[gray!50!magenta, opacity=0.5]
  (0 + \xshift,     0.08073) -- 
  (1*\dx + \xshift, 0.22452) --
  (2*\dx + \xshift, 0.13351) --
  (3*\dx + \xshift, 0.87745) --
  (4*\dx + \xshift, 0.40500) --
  (4*\dx + \xshift, 0.29890) -- 
  (3*\dx + \xshift, 0.87581) -- 
  (2*\dx + \xshift, 0.00000) --
  (1*\dx + \xshift, 0.00013) --
  (\xshift,         0.01837) -- cycle; 
  
  \draw[thick,black]  
  (0 + \xshift,     0.0196) node{x} node[above right]{\scriptsize 0.0196}-- 
  (1*\dx + \xshift, 0.0130) node{x} node[above right]{\scriptsize 0.0130}-- 
  (2*\dx + \xshift, 0.0087) node{x} node[above right]{\scriptsize 0.0087}-- 
  (3*\dx + \xshift, 0.8765) node{x} node[above right]{\scriptsize 0.8765}-- 
  (4*\dx + \xshift, 0.3309) node{x} node[above right]{\scriptsize 0.3309};


\end{tikzpicture}
    \caption{Evaluation of the interdependence for different objectives $F_o$ for the operating points $\mathcal{X}_\alpha$ for $\alpha_o$ provided by Table~\ref{tab:alpha} in black compared to the smallest and largest values provided by Table~\ref{tab:res} in magenta excluding $\mathcal{X}_\text{initial}$.}
    \label{fig:opt}
\end{figure}

The results of this analysis are now applied to derive a reasonable weighting for the cost function.
The weighting parameters are calculated via~(\ref{e:atilde}) using $\mu_{F_o}$ listed in Table~\ref{tab:alpha}.
As displayed in Fig.~\ref{fig:norm:obj}, the objective $\mathcal{X}^*_{\mathcal{B},U}$ performs also well regarding the evaluation on other objectives, since the cyan area in Fig.~\ref{fig:norm:obj} is rather small.
Thus it gets a high rating rating of $\tilde{\alpha}_{\mathcal{B},U} = 10$.
Reactive power injections into the considered and the external grids are weighted equally with $\tilde{\alpha}_{\mathcal{G},Q} = \tilde{\alpha}_{\mathcal{E},Q} = 5$.
A small value for network losses $\tilde{\alpha}_{\text{slack},P} = 1$
is used, since the influence on this objective is rather small, as the standard deviation suggests.
The objective $\mathcal{X}^*_{\mathcal{L},I_\text{S}}$ is omitted, i.~e. $\tilde{\alpha}_{\mathcal{L},I_\text{S}} = 0$, since it has a strong negative effect on the other objectives, as shown by the large orange area in Fig.~\ref{fig:norm:obj}.
The discussed weighting parameters are displayed in Table~\ref{tab:alpha} with the corresponding optimisation results.
\begin{table}[h]
\caption{Educated choice of normalised weighting parameters $\tilde{\alpha}_o$ and corresponding weighting parameters $\alpha_o$ for the objective function with optimisation results. 
}
\centering
\begin{tabularx}{\linewidth}{lXXXXX}
\toprule
$o$
& ${\mathcal{B},U}$
& ${\mathcal{G},Q}$
& ${\mathcal{E},Q}$
& ${\text{slack},P}$
& ${\mathcal{L},I_\text{S}}$
\\
\midrule
$\tilde{\alpha}_o$ & 10 & 5 & 5 & 1 & 0
\\
$\alpha_o$ & 183.39 & 44.765 & 79.695 & 1.1406 & 0.0000
\\
$F_o(\mathcal{X}^*_\alpha)$ & 0.0196 & 0.0130 & 0.0087 & 0.8765 & 0.3309
\\
\bottomrule
\end{tabularx}
\label{tab:alpha}
\end{table}
The optimisation results are also displayed in Fig.~\ref{fig:opt}.
The magenta area in Fig.~\ref{fig:opt} is spanned between the smallest and largest values provided by Table~\ref{tab:res} excluding $\mathcal{X}_\text{initial}$ to only consider feasible operation points.
Fig.~\ref{fig:opt} and Table~\ref{tab:alpha} confirm a good performance regarding all considered objectives.

\begin{figure}[b]
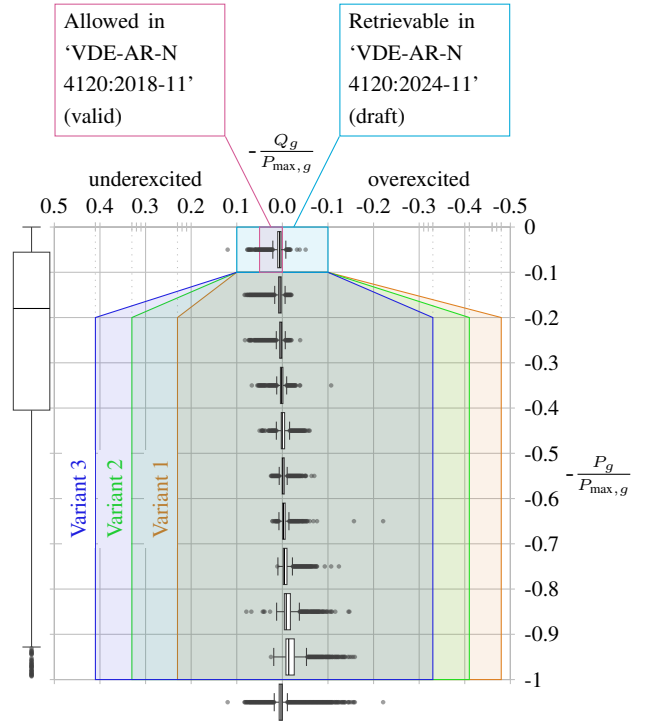

\centering

\caption{Active power dependent evaluation of the 
reactive power feed-ins $Q_g(P_g),g\in \mathcal{G}$ for the operating points $\mathcal{X}_\alpha$ for $\alpha_o,o\in \mathcal{O} $ provided by Table~\ref{tab:alpha} with the $PQ$-variants of~\cite{4120}.}
\label{fig:pq_plot}
\end{figure}

For the remainder of this section, the optimal operation points are investigated for the solution attained by the weighting presented in Table~\ref{tab:res} to classify the results in detail.
The optimal reactive power provisions are displayed together with the corresponding real power injections for all generators in Fig.~\ref{fig:pq_plot}.
A total of 57000 of $PQ$-pairs is analysed.
The abscissa in Fig.~\ref{fig:pq_plot} shows the active power and the ordinate the reactive power in consumer frame convention according to~\cite{4120}.
The SO can choose between 3 variants for grid integration of generators and retrieve reactive power in the coloured regions~\cite{4120}.
An exception to this rule is the magenta region, which corresponds to the current standard `4120:2018-11'~\cite{4120}, within which power plant operators have to operate.
Setpoints aren`t applicable in this region~\cite{4120}.
The corresponding domain of the standard's latest draft, i.~e. `4120:2024-11', is also shown in Fig~\ref{fig:pq_plot}.
According to the draft, reactive power is retrievable within
\begin{equation}
-0.1 P_{\text{max},g}\le    Q_g\le 0.1 P_{\text{max},g}.
\end{equation}

To evaluate the operation points, box-plots are used.
The box-plot on the left shows the distribution of active power feed ins. 
The first quantile for the distribution of ${P_g}/{P_{\text{max},g}}$
is 0.056, the median is 0.1798 and the third quantile is 0.4043.
Thus, during much of the evaluation period, only a small portion of the available active power potential is injected into the grid.

The reactive power provision is analysed in Fig.~\ref{fig:pq_plot} for every real power decile.
Except for a shift towards overexcitedness for high real power injections, the profile is mostly independent of the considered decile.
The quantiles for the distribution of ${Q_g}/{P_{\text{max},g}}$ for the whole data set are -0.0178 for the first quantile, -0.0070 for the median and -0.0030 for the third quantile.
The fences of the distribution are -0.0178 and 0.011 and the percentages of under- and overexcited outliers are 6.47\% and  7.76\%, respectively.

\begin{figure}[ht]
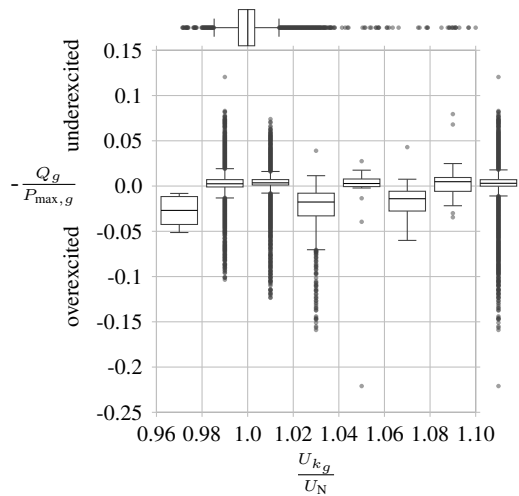

\centering

\caption{Voltage dependent evaluation of the 
reactive power feed-ins $Q_g\left(U_{k_g}\right),g\in \mathcal{G}$ for the operating points $\mathcal{X}_\alpha$ for $\alpha_o,o\in \mathcal{O} $ provided by Table~\ref{tab:alpha}.}
\label{fig:qu_plot}
\end{figure}

The reactive power injections are plotted against the voltage in Fig.~\ref{fig:qu_plot} for generators.
The variants displayed in Fig~\ref{fig:pq_plot} are omitted in Fig.~\ref{fig:qu_plot}, since there are no additional restrictions regarding reactive power injections for voltages between 103~kV (0.9363~pu) and 120~kV (1.0909~pu)~\cite{4120}.
The box-plot of voltage deviation on top of Fig.~\ref{fig:qu_plot} displays the deviation of ${U_{k_g}}/{U_{\text{N}}}$ with fences of 0.985 and 1.0137, a first quantile of 0.9960, a median of 1.0001 and a third quantile of 1.0031.
The percentage of under- and overvoltage outliers is 2.08\% and 4.31\%.
These percentages are reducible with a max-operation applied via (\ref{e:fmax}) on the bus-voltage deviation.
This approach should be considered in future work.
In Fig.~\ref{fig:qu_plot}, the vertical box‐plots show large deviations in the high voltage domain, likely reflecting the small sample size due to the exclusive focus on outliers.

\begin{figure}[ht]
\centering
\begin{tikzpicture}[scale=6]
\draw[gray!50] (0.0, 0.05) -- (0.0, -0.06);
\node[below] at (0.0, -0.05) {\footnotesize 0.96};
\draw[gray!50] (0.1667, 0.05) -- (0.1667, -0.06);
\node[below] at (0.1667, -0.05) {\footnotesize 0.98};
\draw[gray!50] (0.3333, 0.05) -- (0.3333, -0.06);
\node[below] at (0.3333, -0.05) {\footnotesize 1.0};
\draw[gray!50] (0.5, 0.05) -- (0.5, -0.06);
\node[below] at (0.5, -0.05) {\footnotesize 1.02};
\draw[gray!50] (0.6667, 0.05) -- (0.6667, -0.06);
\node[below] at (0.6667, -0.05) {\footnotesize 1.04};
\draw[gray!50] (0.8333, 0.05) -- (0.8333, -0.06);
\node[below] at (0.8333, -0.05) {\footnotesize 1.06};
\draw[gray!50] (1.0, 0.05) -- (1.0, -0.06);
\node[below] at (1.0, -0.05) {\footnotesize 1.08};
\draw[gray] (0,-0.05) rectangle (1,0.05);
\node[below] at (0.5,-0.11) {$\frac{U_{\text{lv},t}}{U_{\text{N-lv},t}}$};
\node[above] at (0.29724999179435047,0.05) {\scriptsize 0.996};
\draw[gray] (0.29724999179435047,0.04) -- (0.29724999179435047,0.06);
\node[above] at (0.09166666666666666,0.05) {\scriptsize 0.971};
\draw[gray] (0.09166666666666666,0.04) -- (0.09166666666666666,0.06);
\node[above] at (0.9833333333333333,0.05) {\scriptsize 1.078};
\draw[gray] (0.9833333333333333,0.04) -- (0.9833333333333333,0.06);
\fill[black!50!gray, opacity = 0.5] (0.175,0.0) circle (0.005);
\fill[black!50!gray, opacity = 0.5] (0.1,0.0) circle (0.005);
\fill[black!50!gray, opacity = 0.5] (0.09166666666666666,0.0) circle (0.005);
\fill[black!50!gray, opacity = 0.5] (0.18333333333333332,0.0) circle (0.005);
\fill[black!50!gray, opacity = 0.5] (0.16666666666666666,0.0) circle (0.005);
\fill[black!50!gray, opacity = 0.5] (0.9833333333333333,0.0) circle (0.005);
\fill[black!50!gray, opacity = 0.5] (0.75,0.0) circle (0.005);
\fill[black!50!gray, opacity = 0.5] (0.48333333333333334,0.0) circle (0.005);
\fill[black!50!gray, opacity = 0.5] (0.8,0.0) circle (0.005);
\draw[black!50!gray] (0.20191544718275056,-0.02) -- (0.20191544718275056,0.02);
\draw[black!50!gray] (0.4283998192116515,-0.02) -- (0.4283998192116515,0.02);
\draw[black!50!gray] (0.20191544718275056,0.0) -- (0.28576475324570855,0.0);
\draw[black!50!gray] (0.3473108156353263,0.0) -- (0.4283998192116515,0.0);
\draw[black!50!gray,fill=white] (0.28576475324570855,-0.04) rectangle (0.3473108156353263,0.04);
\draw[black] (0.29724999179435047,-0.04) -- (0.29724999179435047,0.04);
\draw[gray!50] (0.0, -0.35) -- (0.0, -0.46);
\node[below] at (0.0, -0.45) {\footnotesize -0.7};
\draw[gray!50] (0.125, -0.35) -- (0.125, -0.46);
\node[below] at (0.125, -0.45) {\footnotesize -0.6};
\draw[gray!50] (0.25, -0.35) -- (0.25, -0.46);
\node[below] at (0.25, -0.45) {\footnotesize -0.5};
\draw[gray!50] (0.375, -0.35) -- (0.375, -0.46);
\node[below] at (0.375, -0.45) {\footnotesize -0.4};
\draw[gray!50] (0.5, -0.35) -- (0.5, -0.46);
\node[below] at (0.5, -0.45) {\footnotesize -0.3};
\draw[gray!50] (0.625, -0.35) -- (0.625, -0.46);
\node[below] at (0.625, -0.45) {\footnotesize -0.2};
\draw[gray!50] (0.75, -0.35) -- (0.75, -0.46);
\node[below] at (0.75, -0.45) {\footnotesize -0.1};
\draw[gray!50] (0.875, -0.35) -- (0.875, -0.46);
\node[below] at (0.875, -0.45) {\footnotesize 0.0};
\draw[gray!50] (1.0, -0.35) -- (1.0, -0.46);
\node[below] at (1.0, -0.45) {\footnotesize 0.1};
\draw[gray] (0,-0.45) rectangle (1,-0.35);
\node[below] at (0.5,-0.51) {$\frac{P_{\text{hv},t}}{S_{\text{max},t}}$};
\node[above] at (0.7667763248316831,-0.35) {\scriptsize -0.087};
\draw[gray] (0.7667763248316831,-0.36) -- (0.7667763248316831,-0.34);
\node[above] at (0.051249999999999914,-0.35) {\scriptsize -0.659};
\draw[gray] (0.051249999999999914,-0.36) -- (0.051249999999999914,-0.34);
\node[above] at (0.9941714185015029,-0.35) {\scriptsize 0.095};
\draw[gray] (0.9941714185015029,-0.36) -- (0.9941714185015029,-0.34);
\fill[black!50!gray, opacity = 0.5] (0.08124999999999993,-0.4) circle (0.005);
\fill[black!50!gray, opacity = 0.5] (0.06624999999999992,-0.4) circle (0.005);
\fill[black!50!gray, opacity = 0.5] (0.052499999999999915,-0.4) circle (0.005);
\fill[black!50!gray, opacity = 0.5] (0.07374999999999993,-0.4) circle (0.005);
\fill[black!50!gray, opacity = 0.5] (0.07874999999999993,-0.4) circle (0.005);
\fill[black!50!gray, opacity = 0.5] (0.06749999999999992,-0.4) circle (0.005);
\fill[black!50!gray, opacity = 0.5] (0.06499999999999992,-0.4) circle (0.005);
\fill[black!50!gray, opacity = 0.5] (0.07249999999999993,-0.4) circle (0.005);
\fill[black!50!gray, opacity = 0.5] (0.06374999999999992,-0.4) circle (0.005);
\fill[black!50!gray, opacity = 0.5] (0.05749999999999992,-0.4) circle (0.005);
\fill[black!50!gray, opacity = 0.5] (0.05874999999999992,-0.4) circle (0.005);
\fill[black!50!gray, opacity = 0.5] (0.07499999999999993,-0.4) circle (0.005);
\fill[black!50!gray, opacity = 0.5] (0.051249999999999914,-0.4) circle (0.005);
\fill[black!50!gray, opacity = 0.5] (0.07124999999999992,-0.4) circle (0.005);
\fill[black!50!gray, opacity = 0.5] (0.062499999999999924,-0.4) circle (0.005);
\fill[black!50!gray, opacity = 0.5] (0.08249999999999993,-0.4) circle (0.005);
\draw[black!50!gray] (0.092287821650002,-0.42) -- (0.092287821650002,-0.38);
\draw[black!50!gray] (0.9941714185015029,-0.42) -- (0.9941714185015029,-0.38);
\draw[black!50!gray] (0.092287821650002,-0.4) -- (0.5614405775320452,-0.4);
\draw[black!50!gray] (0.877589543021668,-0.4) -- (0.9941714185015029,-0.4);
\draw[black!50!gray,fill=white] (0.5614405775320452,-0.44) rectangle (0.877589543021668,-0.36);
\draw[black] (0.7667763248316831,-0.44) -- (0.7667763248316831,-0.36);

\end{tikzpicture}
\caption{Transformer analysis regarding voltage and active power for the operating points $\mathcal{X}_\alpha$ for $\alpha_o,o\in \mathcal{O} $ provided by Table~\ref{tab:alpha} with minimum, maximum and median marked.}
\label{fig:tap_pos}
\end{figure}
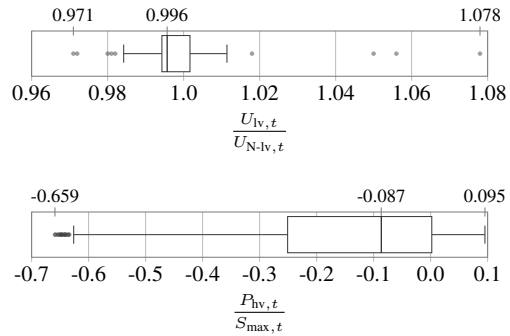

In Fig.~\ref{fig:tap_pos} the characteristic quantities power and voltage are analysed for the grid's transformers.
Fig.~\ref{fig:tap_pos} shows a box-plot for the lower voltage side and underneath a box-plot of active power on the higher voltage side.
The box-plot for the lower voltage side shows a slight, but insignificant undervoltage.
the overvoltage outliers will be reduced in future work by a max-operator applied to the bus voltages as mentioned earlier. 
A voltage dependent evaluation for the high voltage side is irrelevant, since constant voltages of 415~kV (1.092~pu)  and 235~kV (1.068~pu), are present at the nodes of the external grids.

Since the reactive power flow is close to zero, the real power on the transformer's high voltage side is an indicator of the transformer's loading. There are still sufficient reserves, since $|P_{\text{hv},t}|<0.66S_{\text{max},t}$ for the whole data set.

Finally, the tap positions are displayed in Table~\ref{tab:tappos}.
\begin{table}[h]
\caption{Tap positions for the operating points $\mathcal{X}_\alpha$ for $\alpha_o,o\in \mathcal{O} $ provided by Table~\ref{tab:alpha}.
}
\centering
\begin{tabularx}{\linewidth}{lXXXXXX}
\toprule
$\psi_t,\forall t \in\mathcal{T}$ & 6 & 7 & 8 & 9 & 10 & 11
\\
\midrule
Occurrence in \% & 3.70 & 8.44 & 16.33 & 27.66 & 43.04 & 0.70
\\
\bottomrule
\end{tabularx}
\label{tab:tappos}
\end{table}
Only the outliers with an occurrence of one are omitted. 
Table~\ref{tab:tappos} shows that only a small portion of the available tap positions is used.
Currently, the optimiser uses $\psi_{\text{N},t}$ as starting point.
To improve performance in future work the results from \texttt{pandapower}'s power flow should be considered to provide improved initial values for the optimiser.

\section{Conclusion}\label{sec:con}
This paper proposes a novel approach of handling and evaluating multiple objectives for ORPF.
Objectives, considering all operation resources, are investigated.
The interdependence of different objectives is analysed to make an educated choice of weighting parameters for the objective function regarding a specific grid.
The effectiveness of the proposed method is demonstrated on a high voltage grid using 1000 randomly selected study points.
Pareto optimal solutions are derived for all study points with a good performance across all considered objectives.
The subsequent detailed evaluation highlights the importance of the STATCOM domain for reactive power management in the German `Energiewende'.
In future work the max-operator should be applied to penalize outliers regarding individual voltage deviations, while still using the rms-operator to improve the overall voltage profile.
Furthermore, parameters for decentralised control can be derived from the optimisation results using regression technics and the optimiser can be tuned to do economic analysis regarding operations outside of the technical requirements.
Finally, the optimisation should be applied to real networks.

\section{Acknowledgements}\label{sec:ack}
\begin{wrapfigure}[5]{l}{0.17\textwidth}
    \vspace{-0.8cm}
    \begin{center}
        \includegraphics[width=0.17\textwidth]{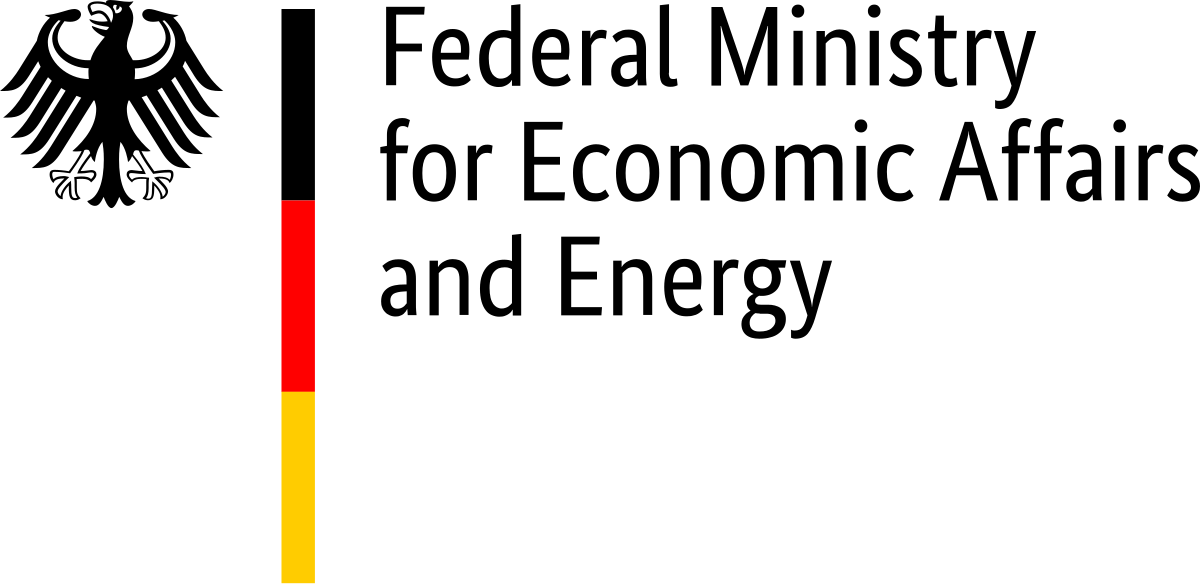}
    \end{center}
\end{wrapfigure}

The investigations concluded for this paper were done in the
context of the research project `Q-REAL' (FKZ 0350061A), which is
supported by the `German Federal Ministry for Economic Affairs and Energy'. 
Authors are solely responsible for the content. 
This paper does not necessarily reflect the consolidated opinion of the project consortium.

\section{Appendix}
Proof of Equations~(\ref{e:Yki}) and~(\ref{e:Ykk}) using the law of conversation of energy, since buses can not store energy:
\begin{align}
    &\;
    \sum_{e\in\mathcal{E}}\mathbbm{1}_{e,k} \underline{S}_e
    +\sum_{g\in\mathcal{G}}\mathbbm{1}_{g,k} \underline{S}_g
    -
    \sum_{m\in\mathcal{M}}\mathbbm{1}_{m,k} \underline{S}_m
    \notag
    \\
    =&\;\underline{U}_k 
    \sum_{i\in\mathcal{B}\setminus k}
    \sum_{{s\in\mathcal{L}\cup\mathcal{T}}}
    \left(
    \mathbbm{1}_{s,(k,i)}
    +
    \mathbbm{1}_{s,(i,k)}
    \right)
    \underline{I}_{s,ki}^*
    \\
    =&\;\underline{U}_k
    \sum_{i\in\mathcal{B}\setminus k}
    \sum_{{s\in\mathcal{L}\cup\mathcal{T}}}
    \Bigg(
    \mathbbm{1}_{s,(k,i)}
    \left(
    \underline{I}_{\text{S},s,ki}^*
    +
    \underline{I}_{\text{P},s,ki}^*
    \right)
    \notag
    \\
    &\;+
    \mathbbm{1}_{s,(i,k)}
    \frac{1}{\underline{n}_s}
    \left(
    \left(\underline{n}_s^*\underline{I}_{\text{S},s,ki}\right)^*
    +
    \left(\underline{n}_s^*\underline{I}_{\text{P},s,ki}\right)^*
    \right)
    \Bigg)
    \\   
    =&\;\underline{U}_k 
    \sum_{i\in\mathcal{B}\setminus k}
    \sum_{{s\in\mathcal{L}\cup\mathcal{T}}}
    \Bigg(
    \mathbbm{1}_{s,(k,i)}
    \left(
    \underline{Y}_{\text{S},s}^*
    \left(
    \underline{U}_k^*-
    \frac{\underline{U}_i^*}{\underline{n}_s^*}
    \right)
    +
    \frac{\underline{Y}_{\text{P},s}^*}2\underline{U}_k^*
    \right)
    \notag
    \\
    &\;+
    \mathbbm{1}_{s,(i,k)}
    \frac{1}{\underline{n}_s}
    \left(
    \underline{Y}_{\text{S},s}^*
    \left(
    \frac{\underline{U}_k^*}{\underline{n}_s^*}-
    \underline{U}_i^*
    \right)
    +
    \frac{\underline{Y}_{\text{P},s}^*}2\frac{\underline{U}_k^*}{\underline{n}_s^*}
    \right)
    \Bigg)
    \\   
    =&\;\underline{U}_k 
    \overbrace{
    \sum_{i\in\mathcal{B}\setminus k}
    \sum_{{s\in\mathcal{L}\cup\mathcal{T}}}
    \left(
    \mathbbm{1}_{s,(k,i)}
    +
    \frac{\mathbbm{1}_{s,(i,k)}}{|\underline{n}_s|^2}
    \right)
    \left(
    \underline{Y}_{\text{S},s}^*
    +\frac{\underline{Y}_{\text{P},s}^*}{2}
    \right)
    }^{=\underline{Y}_{kk}^*}
    \underline{U}_k^*
    \notag
    \\
    &\;    
    +
    \underline{U}_k 
    \sum_{i\in\mathcal{B}\setminus k}
    \underbrace{
    \sum_{{s\in\mathcal{L}\cup\mathcal{T}}}
    -
    \left(
    \frac{\mathbbm{1}_{s,(k,i)}}{\underline{n}_s^*}
    +
    \frac{\mathbbm{1}_{s,(i,k)}}{\underline{n}_s}
    \right)
    \underline{Y}_{\text{S},s}^*
    }_{=\underline{Y}_{ki}^*}
    \underline{U}_i^*
    \\
    =&\;
    \underline{U}_k 
    \sum_{i\in\mathcal{B}}
    \underline{Y}_{ki}^*
    \underline{U}_i^*.
\end{align}
\end{document}